\documentclass{amsart}

\usepackage{
amsmath,amsthm,amsfonts,amscd,graphicx,latexsym,amssymb,eucal,curves}
\usepackage[all]{xy}

\newcommand {\Pfeil} {\longrightarrow}

\newcommand {\ord}{{\rm ord}}

\renewcommand{\iff}{\quad\Leftrightarrow\quad}

\renewcommand{\phi}{\varphi}

\renewcommand{\epsilon}{\varepsilon}

\newcommand{\BS}{{\rm BS}}
\newcommand{\alg}{{\rm alg}}
\newcommand{\an}{{\rm an}}

\newcommand{\orb}{{\rm orb}}

\newcommand{\sharf}{{\,\sharp\,}}

\newcommand{\triv}{{\rm triv}}

\DeclareMathOperator{\Aut}{Aut}
\DeclareMathOperator{\Out}{Out}

\DeclareMathOperator{\Ends}{Ends}
\DeclareMathOperator{\Gal}{Gal}

\DeclareMathOperator{\PGL}{PGL}

\newcommand{\porbS}{{\pi_1^\orb(\mathcal{S},\bar{s})}}

\newtheorem{Def}{Definition}[section]

\newtheorem{Rem}[Def]{Remark}

\newtheorem{Prop}[Def]{Proposition}

\newtheorem*{Satz*}{Satz}

\newtheorem{Lemma}[Def]{Lemma}

\newtheorem{Cor}[Def]{Corollary}

\newtheorem*{Quest*}{Question}

\newtheorem{Thm}{Theorem}

\begin{document}



\title{\bf Riemann Existence Theorems of Mumford Type}

\author{Patrick Erik Bradley}
\date{\today}

\maketitle







\begin{abstract}
Riemann Existence Theorems for Galois covers of Mumford curves by
Mumford curves are stated and proven. As an application,
all finite groups are realised as full automorphism groups
of Mumford curves in characteristic zero.
\end{abstract}

\section{Introduction}

The classical Riemann Existence Theorem (RET) states that
a finite group $G$ can be realised as a deck transformation group of
a $G$-cover of the Riemannian sphere by a compact Riemann
surface, if and only if $G$ is generated by some elements
whose product is one. This amounts to finding
a $\mathbb{C}$-rational point of a Hurwitz space.

\smallskip
If in the situation of $p$-adic geometry we consider
Mumford curves to be the analogue of Riemann surfaces,
the translation ad literam of the RET is no longer valid.
The purpose of this article is to find a correct translation
of the RET in the characteristic zero case. This means that we
relate the non-emptiness of Mumford-Hurwitz spaces
constructed in \cite{BradDiss2002} to certain generating systems
 with the right properties, i.e.\ it is exactly these properties
we will explain in this article.

\smallskip
For this, the notion of Mumford orbifold becomes crucial.
These are Mumford curves covered by Mumford curves,
and to such a cover $\phi$ belongs a certain graph of groups
encoding the ramification behaviour of $\phi$ which in the tree-case
was called {\em $*$-tree} by F.\ Kato who studied them first. 
In the general case, we call it a $*$-graph.
Exploiting its structure 
gives the translation we seek.
We call the result {\em virtual RET}, because of the interpretation of
the stabilisers of vertices in $*$-graphs as 'virtual'
ramification---something that the orbifold itself does not 'see'.

\smallskip
In the special case that all edge groups are
trivial, we recover the so-called {\em Harbater-Mumford components}
of Hurwitz spaces considered by Fried in \cite{Fried1993}.
This leads to another proof of Harbater's result
that every finite group can be realised as the Galois group
of a cover $X\to\mathbb{P}^1_{\mathbb{C}_p}$ of the complex $p$-adic
projective line in characteristic zero by a Mumford curve $X$
\cite{Harbater1987}. He does the job with 'mock covers' which
are totally degenerate covers over the residue field,
i.e.\ they lift to covers by
 Mumford curves.

\smallskip
Until now, it was an open question, if every finite group
can be realised as the full automorphism group of a Mumford curve
in characteristic zero. Our techniques yield a positive answer.

\medskip
The setup for this paper is $p$-adic geometry of characteristic zero,
that is we are working over the field $\mathbb{C}_p$ obtained by completing
an algebraic closure of the field $\mathbb{Q}_p$ of rational $p$-adic
numbers.


\medskip

\section{Hurwitz spaces}

\subsection{Covers of projective curves}    \label{Bertin}

Let $G$ be a finite group of order $n$ and $\phi\colon
X\stackrel{G}{\Pfeil} Y$ a cover of a non-singular projective
irreducible algebraic curve $Y$ of genus $g$ over a field $K$ of
characteristic zero, ramified above $r$ distinct $K$-rational
points $\eta_1,\dots,\eta_r$ with orders $e_1,\dots,e_r$. It is
given by a surjective map of topological groups
$\psi\colon\pi_1^\alg(Y\setminus\{\eta_1,\dots,\eta_r\},y)\to\Aut\phi\cong
G$.

\smallskip
Taking topological generators $\alpha_1,\beta_1,\dots,\alpha_g,\beta_g,\gamma_1,\gamma_r$
of $\pi_1^\alg(Y\setminus\{\eta_1,\dots,\eta_r\},y)$ with the relation
$$
\prod\limits_j[\alpha_j,\beta_j]\prod\limits_i\gamma_i=1,
$$
where $[\alpha,\beta]=\alpha^{-1}\beta^{-1}\alpha\beta$,
we obtain a generating system
$$
a_j:=\psi(\alpha_j),\; b_j:=\psi(\beta_j),\;g_i:=\psi(\gamma_i)
$$
of $G$ with the same relation, called {\em genus $g$ generating system}.
The conjugacy classes $C_i$ of the $g_i$
define the {\em ramification type} ${\bf C}:=(C_1,\dots,C_r)$
of the cover $\phi$.

\medskip
If $K=\mathbb{C}$, then $\gamma_i$ can be viewed as a closed topological path
in $Y\setminus\{\eta_1,\dots,\eta_r\}$ separating $\eta_i$ from the other
branch points, and lifting the path to $X$ gives a deck transformation
$g_i\in \Aut\phi$ which
 permutes
the fibre $\phi^{-1}(y)$  ($y$ is not assumed to be a
branch point). Conjugation means changing the genus $g$
generating system of $G$.

\bigskip
Let $\mathcal{M}_{g,r}$ be the moduli stack of smooth $r$-pointed curves
of genus $g$, and $\mathcal{X}\to\mathcal{M}_{g,r}$ the universal object. For a
finite group $G$ of order $n$ with its regular representation
$G\to S_n$ into the symmetric group, there is, by
\cite[Remark 4.1.3]{WewerDiss}, a smooth algebraic stack $\mathcal{H}^G_\mathcal{X}(G)$ whose
coarse moduli scheme $H_g(G,r)$ is finitely \'etale over the
moduli scheme $M_{g,r}$ of $\mathcal{M}_{g,r}$. This moduli space
parametrises equivalence classes of Galois covers $Y\to X$ of smooth
curves of genus $g$ ramified above $r$ points together with an
isomorphism $\sigma\colon G\to \Aut(Y/X)$. Two pairs
$(Y\stackrel{G}{\longrightarrow}X,\sigma)$ and
$(Y\stackrel{G}{\longrightarrow}X,\sigma)$ are equivalent, if
there is a commutative diagram
$$
\xymatrix{ Y \ar[r]^f \ar[d] & Y' \ar[d] \\
X \ar[r] &X' }
$$
with horizontal isomorphisms, and such that $f\circ\sigma(g)\circ
f^{-1}=\sigma'(g)$ for all $g\in G$.

\smallskip
 Define $\mathcal{H}_\mathcal{X}(G,{\bf C})$ to
be the algebraic substack of $\mathcal{H}^G_\mathcal{X}(G)$ of covers with
ramification type ${\bf C}$ over any irreducible scheme whose function
field is of characteristic zero, and let $H_g(G,{\bf C})$ be its
coarse moduli scheme.  It is a quasi-projective scheme over
$\mathbb{Z}\left[\frac1n\right]$ \cite{Bertin1996}. The space of $\mathbb{C}$-rational
points $H_g(G,{\bf C})(\mathbb{C})$ is non-empty if and only if $G$ has a
genus $g$ generating system. This  is in the genus zero case
nothing but the Riemann Existence Theorem, the staring point of
our observations from Section \ref{mumfloc} on.

\begin{Def} 
The spaces
$${\bf H}_g(G,{\bf C}):=(H_g(G,{\bf C})\otimes_{\mathbb{Z}\left[\frac1n\right]}\mathbb{C}_p)^{\an}\quad
\text{and}\quad
{\bf H}_g(G,r):=(H_g(G,r)\otimes_{\mathbb{Z}\left[\frac1n\right]}\mathbb{C}_p)^\an,
$$
where the superscript $^\an$ means analytification of schemes,
are
called {\em $p$-adic Hurwitz spaces}.
\end{Def}

We will also consider the intermediate Hurwitz spaces ${\bf H}_g(G,{\bf e})$
of covers with fixed list of ramification orders ${\bf e}=(e_1,\dots,e_r)$ which we call the {\em signature} of the covers parametrised by the Hurwitz space.

\subsection{$*$-graphs}
If not stated otherwise, $K$ will always denote a finite extension field of the $p$-adic field $\mathbb{Q}_p$
which is supposed to contain all ramification points of covers.
This is justifiable, as all equivalence classes of covers of curves defined over $\mathbb{C}_p$ 
contain  in fact a representative defined over some finite
extension field of $\mathbb{Q}_p$. 

\medskip
In the non-Archimedean setting one faces the problem that there are not enough
topological coverings, e.g.\ it is a well known fact that all analytic subsets of $\mathbb{P}^1$ or all curves with tree-like
reduction are topologically simply connected. For this reason, more general covers are
often taken into consideration. 

\begin{Def}[De Jong]	\label{etalecover}
An analytic map $f\colon S'\to S$, with $S$ a connected $p$-adic manifold, is an {\em \'etale covering}, if $S$ is covered by open subsets $U$ such that $f^{-1}(U)=\coprod V_j$ and 
$f|_{V_j}$ is a finite \'etale morphism.
\end{Def}

By a {\em $p$-adic manifold} we mean a smooth paracompact strictly $\mathbb{C}_p$-analytic space 
in the sense of \cite{Berkovich1993}. Finite \'etale morphisms are examples of \'etale covers. 
If the $f|_{V_j}$ in Definition \ref{etalecover} are all isomorphisms, then we have {\em topological covers}. These correspond to locally constant sheaves of sets on $S$. In \cite[Lemma 2.2]{deJong1995b} it is shown that
\'etale covers $f\colon S'\to S$ are \'etale and separated morphisms of analytic spaces;  if
$S$ is paracompact then also is $S'$. Thus if $f$ is an \'etale cover of a $p$-adic manifold, then $S'$ is also a $p$-adic manifold. 
According to \cite{Berkovich1999}, such spaces are
locally compact, locally arcwise connected and locally contractible. This means that the
usual theory of universal coverings and topological fundamental groups applies. 

\smallskip
For an \'etale cover $S'\to S$, the \'etale quotient sheaf $S'/R$ on $S$ by an equivalence relation $R\subseteq S'\times_S S'$ is representable by an \'etale cover $S''\to S$, if $R$ is a union of connected components \cite{deJong1995b}. In this case, we call $S''\to S$ a {\em quotient} of $S'\to S$.
 This leads to   
 Andr\'e's notion of {\em tempered}
\'etale covers.

\begin{Def}[Andr\'e]
A {\em tempered} cover $S'\to S$ is a quotient of a composite \'etale cover $T'\to T\to S$,
where $T'\to T$ is a topological covering and $T\to S$ a finite \'etale cover.
\end{Def}

Fixing the ramification orders in covers yields the notion of orbifold
which can quickly be defined as follows: a
(one-dimensional uniformisable) {\em orbifold} 
$\mathcal{S}=(S,(\zeta_i,e_i))$ consists of a $p$-adic manifold $S$ of
dimension one and finitely many points $\zeta_i$ allowing a Galois
cover $\Omega\to S$ (a {\em global chart}) which is tempered
outside $\{\zeta_i\}$ and of ramification orders $e_i$ above
$\zeta_i$.

\begin{Def}	\label{mumfordorbifold}
A {\em Mumford orbifold} is a one-dimensional orbifold
$\mathcal{S}$ having a global chart $\phi\colon X\to S$ with a Mumford curve $X$ (where the projective line and Tate curves are also considered as Mumford curves).
If $\mathcal{S}=(\mathbb{P}^1,(0,e_0),(1,e_1),(\infty,e_\infty))$, then $\mathcal{S}$
will be called a {\em Mumford-Schwarz orbifold}.
\end{Def}

Let $\mathfrak{M}_g$ be the moduli space defined over $\mathbb{C}_p$ of Mumford curves of genus $g$. It is an analytic
subspace of the moduli space $M_g\otimes\mathbb{C}_p$ of all genus $g$ curves.

\begin{Def}	\label{def-MH}
The analytic space defined by the pullback via the cartesian diagram
$$
\xymatrix{
\mathfrak{H}_g(G,{\bf e}) \ar[r]\ar[d] & {\bf H}_g(G,{\bf e})\ar[d] \\
\mathfrak{M}_h \ar[r] & M_h\otimes\mathbb{C}_p
}
$$ 
where the vertical arrows map a cover to the upper curve, and $h$ is given by the Riemann-Hurwitz formula,
is called a {\em Mumford-Hurwitz space}.
\end{Def}
 
Mumford-Hurwitz spaces parametrise 
$G$-covers of Mumford orbifolds of genus $g$ with
fixed lists of ramification orders  ${\bf e}=(e_1,\dots,e_r)$.   
In \cite{BradDiss2002} we  constructed  them 
using Herrlich's non-archimedean Teichm\"uller
spaces from \cite{HerrlichHabil}, and also Andr\'e's orbifold fundamental
group $\porbS$ from \cite{Andreperiod2003}.


\begin{Lemma}	\label{MumfordHurwitzdimension}
The Mumford-Hurwitz space $\mathfrak{H}_g(G,{\bf e})$ is open in ${\bf H}_g(G,{\bf e})$.
\end{Lemma}

\begin{proof}
This follows from the well known fact that $\mathfrak{M}_h$ is open in $M_h\otimes\mathbb{C}_p$ ($h$ as in Definition
\ref{def-MH}), which in turn holds because the fibres of the reduction map into the boundary of the moduli space of curves over $\bar{\mathbb{F}}_p$ are open
in ${\bf H}_g(G,{\bf e})$.
\end{proof}

From Definition \ref{mumfordorbifold}
it follows that if $\mathcal{S}$ is a Mumford orbifold, then $S$ is  a
Mumford curve (e.g.\ \cite[Lemma 4.2]{BradDiss2002}). Let $\Omega\to X$
be the topological universal cover of $X$. The space $\Omega$
is an analytic subdomain of the $p$-adic projective line, and it fits into a
commuting diagram
$$
\xymatrix{
\Omega \ar[r]^{F_g}\ar[dr]_{N}&X\ar[d]^{G}\\
    &\mathcal{S}
}
$$
where $G$ is a finite group and $F_g$ is the free group on $g=\text{genus}(X)$ generators. In fact, $\Omega$ is the complement in $\mathbb{P}^1$ of the limit points
of $N$ acting discontinuously as a subgroup of $\PGL_2(\mathbb{C}_p)$.
According to \cite[4.1.7]{Berkovich1990},
there is an induced action of
the discrete group $N$ on the skeleton $\Delta(\Omega)$, a tree whose ends correspond bijectively
to the set $\mathcal{L}$ of limit points of $N$,
and there is a pure affinoid covering $\mathfrak{U}$ of $S$ such that
the quotient $\Delta(\Omega)/N$ is isomorphic to the
intersection graph of the analytic reduction $\Delta_{\mathfrak{U}}(S)$
 of $S$ with respect to $\mathfrak{U}$. The graph of groups
$\mathfrak{G}(N):=(\Delta_{\mathfrak{U}}(S),N_\bullet)$ contains some information
on the cover $X\stackrel{G}{\Pfeil}\mathcal{S}$: since the topological fundamental
group of $X$  is a free group, all vertex groups of $\mathfrak{G}(N)$ are
contained in $G$. They coincide with the stabilisers
under the action of $G$ on the reduction graph $\Delta_{\phi^{-1}(\mathfrak{U})}(X)$
of $X$. But there is no direct information on the ramification
of the cover $\phi$ itself. As, according to \cite[Satz 5]{Herrlich1980},
the ramification points of $\phi$ are exactly the $F_g$-orbits of fixed points in $\Omega$
of elliptic transformations of $N$,  
Kato remedies this problem by considering
the complement
$\Omega^*$ of the set $\mathcal{L}^*$ which is the union of $\mathcal{L}$ and the set
of fixed points of elements of finite order in $N$.
As $N$ acts also on $\Omega^*\subseteq\Omega$, the graph of groups
$\mathfrak{G}^*:=(\Delta(\Omega^*)/N,N_\bullet)$ is well-defined. Whereas $\mathfrak{G}$ is a finite graph
of groups, $\mathfrak{G}^*$ is an infinite graph with finitely many cusps which correspond bijectively to the branch points of $\phi$ and whose stabilisers coincide with the decomposition groups \cite[Proposition 2.2]{Kato2004}. However, both $\mathfrak{G}$ and $\mathfrak{G}^*$
have fundamental groups isomorphic to $N$.
As all vertex groups are finite, this means by \cite{Khramtsov1991} that
there is a finite number of {\em admissible} edge contractions or {\em slides} 
which turns the finite part $\mathfrak{G}'$  of $\mathfrak{G}^*$ (obtained by contracting the cusps) to $\mathfrak{G}$. Here, we say that a contraction of an edge $e$ in a graph of groups $(\mathfrak{G},N_\bullet)$ is 
{\em admissible}, if $N_e$ is isomorphic to either $N_{o(e)}$ or $N_{t(e)}$, where $o(e)$ is the origin vertex of $e$ and $t(e)$ the terminal vertex of $e$ and $o(e)\neq t(e)$.
A {\em slide} means replacing an embedding $\alpha\colon N_e\to N_v$ (with $v=o(e)$ or $v=t(e)$) by $c_g\circ\alpha$, where $c_g$ is conjugation by $g\in N_v$, but leaving the underlying graph $\mathfrak{G}$ unchanged.


\begin{Def}
$\mathfrak{G}^*$ is called the {\em $*$-graph} of $\phi$.
\end{Def}

If $\mathfrak{G}^*$ is a tree, then we call it a {\em $*$-tree}, of course.
If $X\stackrel{G}{\Pfeil}\mathcal{S}$ is a finite global chart of a Mumford
orbifold with at most three branch points, then the corresponding
$*$-tree is called an {\em elementary $*$-tree}. $G$
is then a finite group of projective linear transformations.



\subsection{The dimension of Mumford-Hurwitz spaces}

By Lemma \ref{MumfordHurwitzdimension}, the dimensionality of the Mumford-Hurwitz space $\mathfrak{H}_g(G,{\bf e})$
is already settled: it is $3g-3+m$, where $m$ is the number of branch points of covers belonging to $\mathfrak{H}_g(G,{\bf e})$. The goal of this subsection is to relate the number of branch points to the structure
of the graph of groups associated to such a cover.

\medskip
In fact, the number of branch points does not depend on the particular choice of 
 the discrete embedding of a finitely generated group into $\PGL_2(K)$:

\begin{Prop}[with H.\ Voskuil] \label{branchindep}
Let $\tau\colon N\hookrightarrow\PGL_2(K)$ be a discrete embedding of a finite tree of finite groups
which are amalgamated along non-trivial finite groups. Then the number $m=\#\Ends(T^*_N)$ does not depend
on the particular choice of $\tau$. Moreover, the stabilisers  of ends of $T^*_N$ (as abstract groups)
do not depend on $\tau$. 
\end{Prop}

\begin{proof}
Consider an end $\mathcal{E}$ of the tree $\mathcal{T}:=\mathcal{T}_{\tau(N)}$ which is contained in the $*$-tree $\mathcal{T}^*_{\tau(N)}$. We first prove that stabilisers of vertices along $\mathcal{E}$ do not form a descending
chain. More precisely, if we take any half line beginning at vertex $v_0$ and whose equivalence class is $\mathcal{E}$,
then the stabilisers $N_0$, $N_1$, \dots of the consecutive vertices $v_0$, $v_1$, \dots  does not descend.
Otherwise, the sequence $N_0\supseteq N_1\supseteq \dots$ would become stationary from a vertex $v_i$ on with smallest group $N_{v_i}=C_\ell$ cyclic. By assumption, $N_{v_i}$ must be non-trivial. 
This means that there is an element $\alpha\in N\setminus\{1\}$  of finite order fixing infinitely many
vertices of $\mathcal{T}$. By \cite[Satz 6]{Herrlich1980}, $\alpha$ commutes with an hyperbolic element $\gamma\in\tau(N)\setminus\{1\}$. But then the quotient graph $\mathcal{T}/N$ contains a loop, which is impossible. 

\smallskip
It follows from the above that if one takes a contraction $\tilde{\mathcal{T}}$ of $\mathcal{T}$ such that for all maximal finite
subgroups $H\subseteq N$ there is a vertex $v$ in $\tilde{\mathcal{T}}$ stabilised by $H$, then to each end of $\tilde{\mathcal{T}}$
there corresponds a unique end of $\mathcal{T}$. Moreover, stabilisers of ends do not change under this correspondence. Since for any discrete embedding $\tau'\colon N\to \PGL_2(K)$ the tree $\mathcal{T}'=\mathcal{T}_{\tau'(N)}$ can be contracted to $\tilde{\mathcal{T}}$, one concludes (using again \cite[Satz 6]{Herrlich1980}) that for any finite-order
element $\alpha\in N\setminus\{1\}$ the assertion
$$
\text{fixed points of $\tau(\alpha)$ are regular} \iff \text{fixed points of $\tau'(\alpha)$ are regular}
$$
holds. This proves  the proposition, as the ends of $T^*_{\tau(N)}$ are the $N$-orbits of 
regular fixed points in $\mathbb{P}^1_K$ of elements of $N\setminus\{1\}$ of finite order.
\end{proof}


\begin{Def} \label{stable}
A graph of groups $(\mathfrak{G},N_\bullet)$ with finite vertex groups is called {\em stable}, if for each vertex $v$
of valency $\le 2$ and each edge $e$ with origin $v$ the inclusion $N_e\to N_v$
maps $N_e$ to a proper subgroup of $N_v$.
\end{Def}







{\sc Amalgamification}

\medskip
Let $N$ be a finitely generated group, and $K$ a complete
non-archimedean valued field. Denote further $ T(N)$ the
space of representations $\tau\colon N\to \PGL_2(K)$ such that
$\tau$ is injective, and $\tau(N)$ acts discontinuously and
does not contain parabolic elements. $T(N)$ is invariant
under action of $\PGL_2(K)$, and the quotient $\bar{T}(N)$ is
called the {\em Teichm\"uller space} for $N$ over $K$. It is
a fine moduli space, and the quotient
$\mathfrak{M}(N):=\bar{T}(N)/\Out(N)$, where the
outer automorphism group $\Out(N)$ acts discontinuously,
parametrises $N$-uniformisable Mumford curves
\cite{HerrlichHabilArtikel}.

\smallskip
Mumford curves uniformised
by groups of a certain type allow some control over the branch locus
of the corresponding covers, as the following will reveal.

\begin{Def}
We say that a group which is a free product of a free group of
finite rank with a free amalgamated product of finitely many
finite groups is {\em of type} {\bf A}{\bf m}. 
\end{Def}

\smallskip
Let $N=G_1*_H G_2$ be an amalgam of finitely generated
groups, $\tau\in T(N)$, $\tau_i:=\tau|_{G_i}$ and
$\tau_H:=\tau|_H$. Denote the corresponding sets of ordinary
points in $\mathbb{P}^1$ as $\Omega$, $\Omega_i$ and $\Omega_H$. We
clearly have $\Omega_i\subseteq\Omega_H$ and
$\Omega\subseteq\Omega_1\cap\Omega_2$, where in general equalities
do not hold. We are interested in the branch loci of the covers
$$
f\colon\Omega\to\Omega/N,\quad
f_i\colon\Omega_i\to\Omega_i/G_i,\quad
f_H\colon\Omega_H\to\Omega_H/H
$$
induced by $\tau$. Denote the branch loci of the above covers with
$br(f)$, $br(f_i)$ and $br(f_H)$.

\begin{Lemma}   \label{branchrestrict}
If $N=G_1*_H G_2$ is of type {\bf A}{\bf m},  then
$$
br(f_H)= br(f_H|_{\Omega})\quad\text{and}\quad
br(f_i)=br(f_i|_{\Omega}).
$$
\end{Lemma}

\begin{proof}
Let $x\in br(f_H)$ resp.\ $\in br(f_i)$. It is the $H$-orbit
resp.\ $G_i$-orbit of a point $z\in\Omega_H$ resp.\ $\in\Omega_i$
fixed by an element $\alpha$ of finite order in $H\setminus\{1\}$
resp.\ in $G_i\setminus\{1\}$. If $x$ were not a branch point of
the restriction to $\Omega$, then there would exist a $z$ as
above, but outside $\Omega$. By \cite[Satz 6]{Herrlich1980}, this
means that there exists a hyperbolic $\gamma\in N$ which
commutes with $\alpha$. But then  a stable graph $\mathfrak{G}$ of groups
with fundamental group $N$ contains an edge which is a loop
and is stabilised by $\alpha$. This implies that $N$ cannot
be of type {\bf A}{\bf m}.
\end{proof}

The calculation of the dimension of the spaces $T(N)$ is reduced by
Herrlich to considering the case of type {\bf A}{\bf m}\ by constructing a finite
\'etale map $T(N)\to T(N')$ with $N'$ of type {\bf A}{\bf m}. 

\begin{Lemma}[Herrlich]	\label{amalgamify}
If $T(N)\neq\emptyset$ then there is a subgroup $N'\subseteq N$ which is of type {\bf A}{\bf m}\
and the canonical map $T(N)\to T(N')$ is finite and unramified.
\end{Lemma}

\begin{proof}
Let $\mathfrak{G}$ be a stable graph of groups with fundamental
group $N$. Following the proof of \cite[Lemma
14]{HerrlichHabilArtikel}, we need to consider only the case,
where the underlying graph $Q$ of $\mathfrak{G}$ is such that for every
 maximal
subtree $P$ the complement  $Q\setminus P$ contains an edge $e$ with
nontrivial stabiliser $N_e$. There are only two possible cases:

\begin{itemize}
\item[(Am 1)] The edge $e$ is a loop, and $N_e=N_{o(e)}$ is cyclic.
In this case, let  $\mathfrak{G}'$ be the graph of groups obtained
from $\mathfrak{G}$ by replacing $N(e)$ and its images under
all embeddings of edge groups by the trivial group. Then $N$
is generated by the  fundamental group $N'$ of
$\mathfrak{G}'$ and an element $\alpha$ commuting with  the
element $\gamma_e$ belonging to the edge $e$. So, any $\tau\in
T(N')$ can be extended in finitely many ways to an element
$\tilde{\tau}\in T(N)$ such that $\tilde{\tau}(\alpha)$ is
non-trivial of finite order. This gives a finite and unramified map
$T(N)\to T(N')$.
\item[(Am 2)] The vertex group $N_{o(e)}=:N_1$ is not cyclic.
In this case, $\mathfrak{G}$ is replaced by a graph of groups
$\mathfrak{G}'$ by giving $e$ a new vertex $v$ as origin, and setting
$N'_v=N_1$. Viewing the fundamental group $N'$ of
$\mathfrak{G}'$ as a subgroup of $N$, we have in the bigger
group the relation $N'_v=\gamma_e N_1\gamma_e^{-1}$,
which means that the induced map $T(N)\to T(N')$
is again finite and unramified.
%
\end{itemize}
\end{proof}

\begin{Def}
The group $N'$ of Lemma \ref{amalgamify} is called an {\em amalgamification} of $N$.
\end{Def}




\begin{Thm}     \label{branchformula}
Let $N$ be a finitely generated group containing a free subgroup of rank $2$, and $\mathfrak{G}$ a stable graph of groups whose fundamental group is 
$N$.
Let further  $\tau\in T(N)$. If $\Omega\subseteq\mathbb{P}^1$ is the
set of ordinary points of $\tau(N)$, then the cover
$\Omega\to\Omega/N$ of the Mumford curve $\Omega/N$ has
exactly $$n=2(C-c)+3(D-d)$$ branch points, where $C$ resp.\ $D$ denotes the number of non-trivial cyclic
resp.\ non-cyclic vertex groups and $c$ resp.\ $d$ the number of non-trivial cyclic resp.\ non-cyclic edge groups of $\mathfrak{G}$.
\end{Thm}

\begin{proof}
We
will
first
show how $n$ does change if we replace $N$ by an amalgamification $N'$.
We use the notations from the proof of Lemma \ref{amalgamify}

\smallskip
In the case (Am 1), $\tau(\alpha)$  
 has the same fixed points as the
hyperbolic $\tau(\gamma)$, which means that these are necessarily limit points of $\tau(N)$ but not of $\tau(N')$.
Therefore the induced covers
$$
f\colon\Omega\to\Omega/N\quad\text{and}\quad
f'\colon\Omega'\to \Omega'/N',
$$
where $\Omega$ resp.\ $\Omega'$ are the sets of ordinary points in
$\mathbb{P}^1$ of $\tilde\tau(N)$ resp.\ $\tau(N')$,
have a difference of 2 in the number of branch points.

\smallskip
In the case (Am 2), $N$ is generated by $N'$, $\gamma_e$ and $N_1$.
As the latter group is conjugated to $N_v'$,  the number of branch points goes up by 3.

\smallskip

Let now $N=G_1*_H G_2$ be of type {\bf A}{\bf m}.  Let $f$
correspond to a representation $\tau \in T(N)$, and $f_i$
resp.\ $f_H$ the covers corresponding to $\tau|_{G_i}$ resp.\
$\tau|_H$. By Lemma \ref{branchrestrict}, the commutative square
$$
\xymatrix{\Omega/H \ar[r] \ar[d] & \Omega/G_1 \ar[d] \\
            \Omega/G_2 \ar[r] & \Omega/N
}
$$
induced by $\tau$ yields a commutative square of maps between the
branch loci
\begin{align} \label{branchsquare}
\xymatrix{ br(f_H) \ar[r] \ar[d] & br(f_1) \ar[d]^{\phi_1} \\
            br(f_2) \ar[r]_{\phi_2} & br(f)}
\end{align}
This yields the inequality
$$
n\ge \# br(f)=:m, \quad\text{i.e.}\quad 3g-3+n\ge3g-3+m,
$$
because $br(f)$ is the union of the images of $\phi_1$ and
$\phi_2$: a branch point of $f$ is the $N$-orbit of  a fixed
point in $\Omega$ of an element $\alpha\in N$ of finite order
which lies in a conjugate of $G_1$ or $G_2$.

\smallskip
 Let now $F_h$ be a finitely generated free normal
subgroup in $N$ of finite index and rank $\ge 2$. The moduli space
$\mathfrak{M}(N,F_h)=\bar{T}(N)/\Out_{F_h}(N)$,
where $\Out_{F_h}(N)$ is the group of $F_h$-invariant outer
automorphisms of $N$, parametrises commuting diagrams
\begin{align*}
 \xymatrix{
 \Omega\ar[r]\ar[dr]&\Omega/F_h\ar[d]\\
 &\Omega/N
 } 
\end{align*} 
\cite{Herrlich1984}.
By Proposition \ref{branchindep} it follows that $\mathfrak{M}(N,F_h)$ embeds
into 
$\mathfrak{H}_g(N/F_h,m)$. 
This Mumford-Hurwitz space has
dimension $3g-3+m$ by Lemma \ref{MumfordHurwitzdimension}, therefore, we have $\dim
\mathfrak{M}(N,F_h)=3g-3+n\le 3g-3+m$, in other words: $n\le m$.
\end{proof}




\begin{Cor} \label{pushout}
If $G=G_1*_H G_2$ is of type {\bf A}{\bf m}, then
$$
br(f)=br(f_1)\sharf_{br(f_H)}br(f_2),
$$
where $\sharf$ denotes the pushout (here: of sets).
\end{Cor}

\begin{proof}
In the proof of Theorem \ref{branchformula} we have seen that
$br(f)$ is the union of the images of $\phi_1$ and $\phi_2$ in
Diagram (\ref{branchsquare}), so there is a surjective map from
the pushout to $br(f)$.  The equality of the finite cardinalities
of $br(f)$ and the pushout show that the two sets are equal.
\end{proof}

\section{The Mumford curve locus}   \label{mumfloc}
\subsection{The Harbater-Mumford-components}

In \cite{Harbater1987}, Harbater  realises all finite groups $G$ as quotients of
$\Gal(\bar{\mathbb{Q}}_p(T)/\mathbb{Q}_p(T))$ by pasting cyclic covers
$$
\mathbb{P}^l\to\mathbb{P}^l,\;z\mapsto z^n
$$
in a clever way. He thus gets  global
orbifold charts $\phi\colon X\stackrel{/G}{\longrightarrow}\mathcal{S}$
for rational orbifolds
$\mathcal{S}=(\mathbb{P}^1,(\zeta_1,e_1),\dots,(\zeta_{2r},e_{2r}))$ with $e_i=e_{i+r}$
for $i=i,\dots,r$ and a Mumford curve $Y$.
He does this with `mock' covers of $\mathbb{P}^l_{\mathbb{F}_p}$ lifting to totally
degenerated curves ramified above $2r$ points of the projective line,
if the group $G$ is generated by $r$ elements.
The translation of Harbater's proof in \cite[Section 6.4]{BradDiss2002} to $*$-trees
shows that the marked points
have to be separable into pairs $\zeta_i,\zeta_{i+r}$ by a pure affinoid
covering of $\mathbb{P}^l$. This means that from each vertex in the $*$-tree
for the chart $\phi$ either two cusps emanate or none at all.
In this special case, all vertex stabilisers are cyclic, and
all edge groups are trivial.

\smallskip
The ramification data ${\bf C}=(C_1,\dots C_{2r})$ for those covers
with conjugation classes of the deck group $G$
are
required to have a representative of the form
$$
(g_1,\dots,g_r,g_1^{-1},\dots,g_r^{-1}).
$$
This motivates M.~Fried to call it an {\em Harbater-Mumford-representative}
of ${\bf C}$.

\smallskip
In \cite[3.21]{Fried1993} M.~Fried proves that  the
following  condition on $G$ and ${\bf C}$
with a Harbater-Mumford-representative forces the subspace
$\mathcal{H}_{0}(G,{\bf C})$ of the Hurwitz space $\mathcal{H}_{0}(G,2r)$ with
ramification type ${\bf C}$ to be a component (we shall call it
a {\em Harbater-Mumford-component}):

\begin{itemize}
\item[HM] If any pair of inverse conjugation classes is removed from ${\bf C}$,
then the remaining $2r-2$ classes generate the group $G$.
\end{itemize}


\begin{Prop}    \label{HM}
Let $G$ be any finite group. Then Harbater-\-Mumford-\-components of the Hurwitz space
${\bf H}_{0}(G,2r)$ contain covers of rational
Mumford orbifolds by Mumford curves if and only if
$G$ is a quotient of a free tree product of cyclic groups.
\end{Prop}

\begin{proof}
If a Harbater-Mumford component contains covers of Mumford orbifolds by Mumford curves, then
$G$ is a quotient of the fundamental group of a $*$-tree with trivial edge groups, as we have seen in the
beginning of this subsection.

\medskip
Let $G=\Pi/\Sigma$ be a quotient of a free tree product $\Pi$ of cyclic groups.
The group $\Pi$ is the fundamental group of a finite tree with cyclic vertex groups and trivial edge groups
which can be
made into a discontinuously embeddable tree by subdividing segments of the type as in Figure
\ref{embedsegment}.

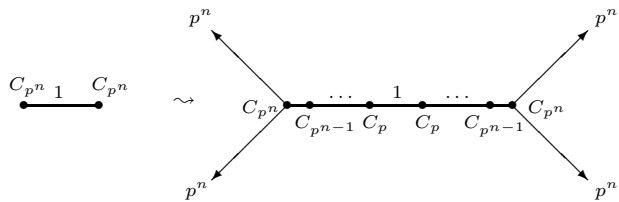
\begin{figure}[h]
\begin{center}
\setlength{\unitlength}{1cm}
\begin{picture}(12,4)
\scriptsize
\put(1,2){\circle*{.1}}
\put(1,2){\line(1,0){1}}
\put(2,2){\circle*{.1}}
\put(.8,2.2){$C_{p^n}$}
\put(1.9,2.2){$C_{p^n}$}
\put(1.4,2.1){$1$}

\put(3,2){$\leadsto$}

\put(4.5,2){\circle*{.1}}
\put(4.5,2){\line(1,0){3}}
\put(7.5,2){\circle*{.1}}
\put(3.9,1.9){$C_{p^n}$}
\put(7.7,1.9){$C_{p^n}$}
\put(5.9,2.1){$1$}

\put(4.5,2){\vector(-1,1){1}}
\put(4.5,2){\vector(-1,-1){1}}
\put(7.5,2){\vector(1,1){1}}
\put(7.5,2){\vector(1,-1){1}}

\put(5.6,2){\circle*{.1}}
\put(6.3,2){\circle*{.1}}
\put(5.5,1.7){$C_p$}
\put(6.2,1.7){$C_p$}

\put(4.8,2){\circle*{.1}}
\put(7.2,2){\circle*{.1}}
\put(4.6,1.7){$C_{p^{n-1}}$}
\put(6.85,1.7){$C_{p^{n-1}}$}

\put(5.05,2.1){\dots}
\put(6.6,2.1){\dots}

\put(3.2,3.1){$p^n$}
\put(3.15,.8){$p^n$}
\put(8.6,3.1){$p^n$}
\put(8.6,.8){$p^n$}
\end{picture}
\end{center}
\caption{Subdivision of a segment gives a discontinuously embeddable tree with isomorphic fundamental group.}
\label{embedsegment}
\end{figure}

The non-trivial edge groups are the intersections of the stabilisers of the edges' extremities.
The tree on the right of Figure \ref{embedsegment}
is realizable because it is the connected sum of two copies of the $*$-tree
for $C_{p^n}$ with an allowed segment: Indeed, the segment with the trivial edge group
belongs to Herrlich's list \cite{Herrlich1982} of segments of groups  embeddable into the
Bruhat-Tits tree. The geodesics left and right are the $*$-trees of $C_{p^n}$,
and the paths joining the segment to the other two trees are nothing but a part of
the fixed point ``strip'' of a $p$-group in the Bruhat-Tits tree.

\smallskip
More generally, making explicit what Kato proposes in \cite[8.1]{Kato2004}, a segment whose
fundamental group is the free product
$C_e*C_{e'}$ with $e=m p^r$, $e'=m' p^{r'}$ and $\gcd(m,p)=\gcd(m',p)=1$ can be made realisable
by subdividing in the way as in Figure \ref{generalsegment}.

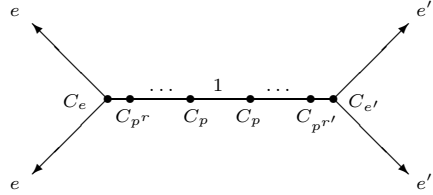
\begin{figure}[h]
\begin{center}
\setlength{\unitlength}{1cm}
\begin{picture}(9,4)
\scriptsize
\put(2.5,2){\circle*{.1}}
\put(2.5,2){\line(1,0){3}}
\put(5.5,2){\circle*{.1}}
\put(1.9,1.9){$C_{e}$}
\put(5.7,1.9){$C_{e'}$}
\put(3.9,2.1){$1$}

\put(2.5,2){\vector(-1,1){1}}
\put(2.5,2){\vector(-1,-1){1}}
\put(5.5,2){\vector(1,1){1}}
\put(5.5,2){\vector(1,-1){1}}

\put(4.4,2){\circle*{.1}}
\put(3.6,2){\circle*{.1}}
\put(3.5,1.7){$C_p$}
\put(4.2,1.7){$C_p$}

\put(2.8,2){\circle*{.1}}
\put(5.2,2){\circle*{.1}}
\put(2.6,1.7){$C_{p^{r}}$}
\put(5,1.7){$C_{p^{r'}}$}

\put(3.05,2.1){\dots}
\put(4.6,2.1){\dots}

\put(1.2,3.1){$e$}
\put(1.2,.8){$e$}
\put(6.6,3.1){$e'$}
\put(6.6,.8){$e'$}
\end{picture}
\end{center}
\caption{A discontinuous embedding of an amalgam of cyclic groups.}
\label{generalsegment}
\end{figure}

The reason is that the two decomposition groups contain $p$-groups with thick fixed point strips.

\smallskip
All other segments with vertex groups of order prime to $p$  are realisable \cite{Herrlich1982}.
Thus, we have made the tree into a $*$-tree whose fundamental group embeds discontinuously
into $\PGL_2(\mathbb{C}_p)$. Let $\Omega\subseteq\mathbb{P}^1$ be its domain of regularity. It is also the domain of regularity for the normal subgroup $\Sigma$ of finite index. Then we
have a commuting diagram
$$
\xymatrix{
\Omega\ar[r]^\Sigma\ar[dr]_\Pi&X\ar[d]^G\\
 & \mathbb{P}^1
}
$$  
where $X$ is a Mumford curve.
\end{proof}

\subsection{$*$-trees with trivial edge groups}

First, we consider $*$-trees $T$  having only trivial edge groups.
As the Bass-Serre fundamental group $\pi_1^\BS(T)$ of a tree of groups
is generated by the vertex groups, this is also the case for
all quotients of $\pi_1^\BS(T)$. In particular, our Galois group $G$
is generated by generators of the vertex groups, which we interpret as
local monodromies in the following way:

\smallskip
$T$ is obtained by pasting together local data:  elementary
$*$-trees are glued by connecting them with edges  at appropriate
vertices. Let now $T'$ be a stable tree of groups obtained from
$T$ by contracting edges whose stabiliser equals one of its extremities'
groups. Now, by Belyi's Theorem the Mumford-Schwarz-orbifolds corresponding
to those elementary $*$-trees are defined over $\bar{\mathbb{Q}}$.
This gives us local (complex) monodromies $\gamma_0,\gamma_\infty$
resp.~$\gamma_0,\gamma_1,\gamma_\infty$ around the cusps
$0,\infty$ resp.~$0,1,\infty$, obeying the relations
$$
\gamma_0\gamma_1=1\quad\text{resp.~}\quad \gamma_0\gamma_1\gamma_\infty=1.
$$

\smallskip
This motivates the following

\begin{Def}
Let ${\bf C}=(C_1,\dots,C_r)$ be a generating system of conjugation classes
of $G$ having the property that there exist representatives $\gamma_i\in C_i$
and a partition of the set $\{1,\dots, r\}$ into pairs and triples
in such a way that for each pair $i,j$ resp.\ each triple $i,j,k$ the relation
$$
\gamma_i\gamma_j=1\quad\text{resp.}\quad\gamma_i\gamma_j\gamma_k=1
$$
holds and all $\langle \gamma_i,\gamma_j\gamma_k\rangle$ are isomorphic to subgroups of $\PGL_2(\mathbb{C}_p)$.
This ramification datum ${\bf C}$ is called
{\em of Mumford type}.
\end{Def}

We may now formulate a special case of a Riemann Existence Theorem for
Mumford curves:

\begin{Thm}[RET of Mumford type] \label{RETmumford}
Let $G$ be a finite group.
Then the Mumford-Hurwitz space $\mathfrak{H}_0(G,{\bf C})$ is non-empty if ${\bf C}$ is of Mumford type.
\end{Thm}

\begin{proof}
Let ${\bf C}$ be of Mumford type. Taking any connected sum of $*$-trees for the groups
generated by the chosen pairs and triples of $\gamma_i\in C_i$ as described above,
we get a tree $T$ of groups which can be made discontinuously embeddable 
by subdividing similarly as in the
proof of Proposition \ref{HM}. If we denote by $N$ the fundamental group of $T$, then $G$ is
obviously a quotient of $N$, and the Galois chart $\phi\colon X\to\mathcal{S}$ is constructed
in an analogous way as the HM-cover from the proof of Prosposition \ref{HM}. 
\end{proof}

\subsection{Virtual ramification in $*$-trees}

The motivation for this subsection is the example \cite[III.6.4.6]{Andreperiod2003}
of a Mumford-Schwarz-orbifold $\mathcal{S}$ whose orbifold fundamental group
$\porbS$
is not topologically generated by its local monodromies $\gamma_0,\gamma_1$ around the
points $0,1$.
Neither does there exist 
a local monodromy $\gamma_\infty$ around $\infty$ such that $\gamma_0,\gamma_1,\gamma_\infty$ obey the relation 
$\gamma_0\gamma_1\gamma_\infty=1$. The reason for this is that already
a quotient of $\porbS$  isomorphic to a fundamental  group of a $*$-tree with three cusps, does not allow this.
A discrete quotient of $\porbS$ is called a {\em $p$-triangle group (of Mumford type)}.  

\bigskip
In \cite{Kato99} there is a first attempt towards classifying all discrete $p$-adic triangle groups.
Refining the methods used there yields a complete classification of all discrete subgroups of 
$\PGL_2(K)$ for  finite extensions $K$ of $\mathbb{Q}_p$ \cite{BKVdiscretePGL2}.
In what follows we will make use of this classification result.



\begin{Prop}		\label{Harmclassify}
Every indecomposable discrete subgroup $N$ of $\PGL_2(K)$ generated by elements of finite order is an amalgamation which corresponds to a tree of groups ${\bf T}(N)$ with the following properties:
\begin{enumerate}
\item The vertex groups $N_v$ are finite non-cyclic groups or $p$-adic triangle groups.
\item The edges $e$ correspond to  maximal cyclic subgroups of $N_{o(e)}$ resp.\ $N_{t(e)}$.
\item For each branch point $x$ of $N_v$ there is at most one edge $e$ with $o(e)=v$, and
$N_e$ is the stabiliser of $x$.
\end{enumerate} 
\end{Prop}


\begin{proof}[Sketch of proof]
(Cf.\ \cite{BKVdiscretePGL2}). Suppose, that there is an indecomposable discrete subgroup $N$ of $\PGL_2(K)$ not satisfying the properties of the proposition. It is a finite amalgam of $m_N$ maximal finite subgroups. Choose $N$ such that $m_N$ is minimal under all counterexamples.

\smallskip
We claim that no pair of finite groups occurring in the amalgam is amalgamated along a cyclic group unless one of the finite groups is itself cyclic.

\smallskip
One sees that there are no two maximal finite subgroups in this amalgam which are amalgamated along a 
common maximal cyclic subgroup, otherwise one would easily obtain a counterexample $N'$ with $m_{N'}<m_N$, a contradiction.

\smallskip
It remains the case that two maximal finite non-cyclic subgroups $G_1$, $G_2$ in the amalgam are amalgamated along a common cyclic subgroup $G_0$ which is not maximal in one of them. One can then see that
$G_1*_{G_0} G_2$ is an amalgam of a finite group and a triangle group along a common maximal cyclic subgroup. 
Again it follows that $m_N$ is not minimal, a contradiction.

\smallskip
As a consequence, $N$ is isomorphic to an amalgam of finite non-cyclic groups along non-cyclic subgroups
amalgamated also with cyclic groups along cyclic groups. Such a group is necessarily a triangle group.   
\end{proof}

In \cite{BKVdiscretePGL2} there is also a complete classification of all discrete $p$-adic triangle groups.

\begin{Def}
The amalgam and the tree of groups from Proposition \ref{Harmclassify} are called {\em regular}.
\end{Def}

\begin{Cor}
The number of branch points for an indecomposable discrete subgroup $N$ of $\PGL_2(K)$ is the number of vertices in  ${\bf T}(N)$ plus two.
\end{Cor}

\begin{proof}
Note that all vertices of the regular tree give a contribution of three branch points. The statement now follows easily from Theorem \ref{branchformula}  together with the fact   that $N$ is the fundamental group of a tree of groups.  
\end{proof}





\begin{Prop}[F.\ Kato, F.\ Herrlich] \label{notriangles}
 There are no infinite discrete $p$-adic triangle groups, if $p>5$. For $p\le 5$ there are infinitely many.
  \end{Prop}

\begin{proof}
The first statement is a result of the classification in \cite{BKVdiscretePGL2}. The second statement follows from the classification of discontinuously embeddable segments of finite groups \cite{Herrlich1982}.
\end{proof}

\begin{Def}
A system of conjugation classes ${\bf C}=(C_0,C_1,C_\infty)$ of a finite group $G$ is called {\em of Mumford-Schwarz type}, if there is a $G$-cover $X\stackrel{G}{\longrightarrow}\mathbb{P}^1$ 
with ramification data ${\bf C}$ and with a Mumford curve $X$.
\end{Def}

\bigskip
{\sc Regular amalgams}

\medskip
Let ${\bf T}$ be a regular tree of groups for a regular amalgam $\Gamma$. For any subtree of groups $T_0\subseteq{\bf T}$ denote $c_{\bf T}(T_0)$ the number branch points of the cover $\varpi_\Gamma\colon\Omega\to\mathbb{P}^1$
which are
contributed by  the subtree $T_0$. 

\begin{Lemma}	\label{localcontribution}
Let a regular tree ${\bf T}$ be spanned by disjoint subtrees $T_1$ and $T_2$ through an edge 
$e$. Then
$$
c_{\bf T}(T_1)=c_{T'}(T_1)
$$
for any subtree $T'$ of ${\bf T}$ containing $T_1$ and $e$.
\end{Lemma}

\begin{proof}
As regular amalgams are of type {\bf Am},  Corollary \ref{pushout} says that the set of branch points corresponding to $T'$ is the pushout of the branch loci of $T_1$ and $T_2\cap T'$ along the two branch points corresponding to $e$. Hence, the contribution from $T_1$ does not depend on the complementary tree $T_2\cap T'$.
\end{proof}

\begin{Lemma}
Let ${\bf T}$ be a regular tree of groups, and $v$ a vertex of ${\bf T}$ with non-cyclic stabiliser $N_v$. Then 
$$
c_{\bf T}(v)=3-\ord(v).
$$
In particular, the order of any vertex in ${\bf T}$ is at most three.
\end{Lemma}

\begin{proof}
Let $\ord(v)=m$, which is a non-negative number. The vertex $v$ is connected by $m$ edges to some components $T_1,\dots,T_m$.  We have
%
%
$$
c_{\bf T}({\bf T})=c_{\bf T}(v)+\sum\limits_{i=1}^m c_{\bf T}(T_i)=c_v(v)-2m+\sum\limits_{i=1}^m c_{T_i}(T_i).
$$
The equation on the right holds true because of Theorem \ref{branchformula}.
Now, if $T^i$ is the tree spanned by $T_i$ and $v$, then
$$
c_{T_i}(T_i)=c_{T^i}(T_i)+1,
$$
because 
$$
c_{T^i}(T_i)+c_{T^i}(v)=c_{T^i}(T^i)= c_{T_i}(T_i) - 2 +c_v(v)=c_{T_i}(T_i)+1,
$$
and $c_{T^i}(v)=2$. The latter can be seen by first using Lemma \ref{localcontribution} in order to see that 
$c_{T^i}(v)=c_{[w,v]}(v)$, where $w$ is the vertex in ${\bf T}$ nearest to $v$, and then verifying $c_{[w,v]}(v)=2$
which follows from the fact that the geodesic in the Bruhat-Tits tree between the two fixed points of the edge stabiliser gets folded under the action of $N_v$.
Since, by Lemma \ref{localcontribution}, $c_{T^i}(T_i)=c_{\bf T}(T_i)$,
we have in fact
$$
c_{\bf T}({\bf T})=3-m+\sum\limits_{i=1}^m c_{\bf T}(T_i),
$$
from which we conclude that 
$$
c_{\bf T}(v)=3-m.
$$
This number is non-negative, if and only if $m\le 3$.
\end{proof}

\begin{Def}
A Galois cover $f\colon X\stackrel{G}{\longrightarrow}\mathbb{P}^1$ with a Mumford curve $X$ is called {\em regular},
if the fundamental group $N$ of its $*$-tree is indecomposable. In this case  the induced cover 
$\Omega\stackrel{N}{\longrightarrow}\mathbb{P}^1$ factoring over $f$  is also called {\em regular}.
\end{Def}

\begin{Thm}[Virtual RET] \label{virtualRET}
Let $G$ be a finite group and ${\bf C}=(C_1,\dots,C_r)$ some conjugation classes in $G$
with $\gamma_i\in C_i$.  Then there is a Galois cover $X\stackrel{G}{\longrightarrow}\mathbb{P}^1$ with a Mumford curve $X$ and ramification type ${\bf C}=(C_1,\dots,C_r)$, if and only if there is a partition
$$
 \{1,\dots,r\}=\bigcup\limits_j I_j
$$ 
into sets $I_j$ of cardinality at least two such that for each $I:=I_j$ the tuple $C_I=(C_i)_{i\in I}$
gives rise to a regular cover $\Omega_I\stackrel{\Gamma_I}{\longrightarrow}\mathbb{P}^1$.
\end{Thm}

\begin{proof}
The proof of Theorem \ref{RETmumford} carries over to this more general situtation.
\end{proof}

The structure of regular trees now yields:

\begin{Thm} 	\label{RETregular}
If $X\stackrel{G}{\longrightarrow}\mathbb{P}^1$ is a non-cyclic regular  cover with ramification type
${\bf C}=(C_1,\dots,C_r)$ and $\gamma_i\in C_i$, then there is a  system
${\bf g}=(\gamma_1,\dots,\gamma_r,\gamma_{r+1},\dots,\gamma_{3s})$ of elements in $G$ and a partition
of $\{1,\dots,3s\}$ into triples $I=\{i_0,i_1,i_\infty\}$ such that 
$$
C_I=(C(\gamma_{i_0}),C(\gamma_{i_1}),C(\gamma_{i_\infty}))
$$ 
is of Mumford-Schwarz type
for some $p$-adic triangle group $\Delta_I\stackrel{\phi_I}{\longrightarrow}G$ mapping to $G$.
The conjugation classes in $C_I$ are to be taken in the image of $\phi_I$.
\end{Thm}

\begin{proof}
Let the regular cover be given as
$$
\xymatrix{
\Omega \ar[dr]_N \ar[r]^F& X \ar[d]^G\\
&\mathbb{P}^1
}
$$
with $F$ a finitely generated free group,
and let ${\bf T}$ be the corresponding regular tree. For each vertex $v$ of ${\bf T}$ we have a $p$-adic triangle group
$\Delta_v$ and $c_{\bf T}(v)$ cusps with  generators $\gamma\in G$ of their stabilisers. This locally defines a cover
$\phi_v\colon X_v\stackrel{\Delta_v}{\Pfeil}\mathbb{P}^1$ with $3-c_{\bf T}(v)$ more branch points. Setting 
$F_v:=F\cap \Delta_v$, a finitely generated free group, we obtain a commuting diagram with exact rows
$$
\xymatrix{
1\ar[r]& F_v\ar[r]\ar[d] & \Delta_v \ar[r]\ar[d] & G_v \ar[r]\ar[d]& 1\\
1 \ar[r]&F\ar[r]&N\ar[r]&G\ar[r]&1
}
$$
As the left and the middle vertical arrows are injective, so is also the right vertical map.
Hence, the decomposition groups of $\phi_v$ embed into $G$, and we can complete the $c_{\bf T}(v)$
conjugacy classes of chosen elements of $G$ with further elements of $G$ to a system of Mumford-Schwarz type 
for the local cover $\phi_v$. Doing this for all vertices of ${\bf T}$ yields $s$ Mumford-Schwarz type systems of conjugacy classes in $G$ or, in other words, the required ${\bf g}=(g_1,\dots,g_{3s})$ containing all the 
$\gamma_i$ ($i=1,\dots,r$) together with a partition of $\{1,\dots,r\}$ into triples as in the theorem.   
\end{proof}

\begin{Def}
The system of conjugacy classes ${\bf C}$ of $G$ as in Theorem \ref{RETregular} is called {\em virtually of Mumford type}.
\end{Def}

The reason for calling Theorem \ref{virtualRET} virtual RET lies in the structure of regular trees, which conceal
local monodromies in their vertex groups. This can be seen most easily when there are no triangle
groups present.

\begin{Cor} \label{virtualmonodromy}
If, in the situation of Theorem \ref{RETregular},  all triangle groups $\Delta_I$ 
are finite,
then  the system ${\bf g}$ generates $G$ and  each triple $I=\{i_0,i_1,i_\infty\}$ ($I=I_j$ from a partition as in Theorem \ref{RETregular}) can be chosen such that
the relation 
$$
\prod\limits_{i\in I}\gamma_i =1 
$$ 
holds for some order on $I$.
\end{Cor}

\begin{proof}
Let $N$ be the fundamental group of the tree of groups associated with the regular cover.
As all the triangle groups $\Delta_I$ are finite, the relations 
$$
\prod\limits_{i\in I}\delta_i=1
$$ 
are feasible with elements $\delta_i\in \Delta_I$, hence in $G$. 
As $N$ is generated by all vertex groups of ${\bf T}$, it follows that the set of all $\gamma_i=\phi_I(\delta_i)$ generates $G$.
\end{proof}

\begin{Def}
The $\gamma_i$ from Corollary \ref{virtualmonodromy} are called {\em virtual monodromies}.
\end{Def}

\begin{Rem}
{\rm
As remarked by Y.\ Andr\'e in \cite[Remark III.6.4.7]{Andreperiod2003}, infinite $p$-adic triangle groups
of Mumford type are not generated by their local monodromies (i.e.\ generators of the decomposition groups of the branch points), there exist finite quotients with the same property.
This is why Corollary \ref{virtualmonodromy} does not hold in the  presence of infinite vertex groups.
In the general  case that ${\bf T}$ has infinite triangle groups as vertex groups, an exact description of virtual local monodromies  as generating systems of Galois groups is yet to be found. }
\end{Rem}

\subsection{From trees to graphs}\label{trees-to-graphs}

In order to generalise the preceding to covers
$\phi\colon X\stackrel{G}{\Pfeil}\mathcal{S}$
with an arbitrary
Mumford orbifold $\mathcal{S}$, it is convenient to use the following
setup: let $\mathfrak{G}$ be the $*$-graph corresponding to the global chart
$\phi$ of $\mathcal{S}$, where we assume $X$ to be a Mumford curve of genus $g$.
Further, let $N:=\pi_1^\BS(\mathfrak{G})$ be the Bass-Serre fundamental group
of the graph of groups $\mathfrak{G}$.
Assuming that $N$ is of type {\bf Am},
  there is a commuting
 diagram
with exact rows and surjective vertical arrows
\begin{align}	\label{splitdiagram}
\xymatrix{
1 \ar[r] & \langle N_\bullet\rangle_N \ar[d]\ar[r] &N\ar[r]\ar[d]
    & F_g \ar[r]\ar[d] & 1\\
1 \ar[r]& \langle N_\bullet\rangle_G \ar[r]&G\ar[r]&F\ar[r]&1
}
\end{align}
where $\langle N_\bullet\rangle_N$ is the normal subgroup of $N$ generated by
the vertex and edge groups of $\mathfrak{G}$. It is
isomorphic to the fundamental group of a $*$-tree obtained
by deleting $g$ appropriate edges of $\mathfrak{G}$.
One sees readily that the top row of (\ref{splitdiagram}) splits.

\begin{Def}
If for a Mumford curve $X$ of genus $g$, the group $N$ in the diagram
\begin{align}	\label{Mumforddiagram}
\xymatrix{
\Omega\ar[r]^{F_g} \ar[dr]_N & X \ar[d]^G \\
&Y
}
\end{align}
is of type {\bf Am}, then we say that the Galois cover $X\stackrel{G}{\longrightarrow}Y$ is also of type {\bf Am}.
\end{Def}

Here is now a first generalised virtual Riemann Existence Theorem:

\begin{Thm} \label{general}
The Mumford-Hurwitz space $\mathfrak{H}_g(G,{\bf C})$ contains a cover of type {\bf Am}, if and only if
$G$ contains a normal subgroup $H$ containing representatives of conjugacy classes from ${\bf C}$ such that ${\bf C}$,
viewed as a ramification type for $H$, is virtually
of Mumford type, $G/H$ is a quotient of $F_g$ and the exact sequence
$$
1\to H\to G\to G/H\to 1
$$
splits.
\end{Thm}

\begin{proof} 
$\Rightarrow$.
Let a cover as in (\ref{Mumforddiagram}) be given with $N$ of type {\bf Am}. Then we obtain a commutative diagram with exact rows and columns
\begin{align} \label{nineexact}
\xymatrix{
&1\ar[d]&1\ar[d]&1\ar[d]&\\
1\ar[r]&F_b\ar[d]\ar[r]&F_h\ar[d]\ar[r]&F_c\ar[d]\ar[r]&1\\
1\ar[r]&\langle N_\bullet\rangle_N \ar[d]\ar[r]&N\ar[d]\ar[r]&F_g\ar[d]\ar[r]&1\\
1\ar[r]&H\ar[d]\ar[r]&G\ar[d]\ar[r]&F\ar[d]\ar[r]&1\\
&1&1&1&\\
}
\end{align}
whose lowest row is the second row of (\ref{splitdiagram}). We still need to show that the lowest row of (\ref{nineexact}) splits. As the middle row splits, the images of $N_0:=\langle N_\bullet\rangle_N$ and $F_g$ 
under the middle vertical arrow to $G$ intersect trivially. So, one sees that there is a map $F\to G$ making the lowest row split.

\smallskip
$\Leftarrow$. 
We construct the cover of type {\bf Am} by showing that there is a diagram (\ref{nineexact}).
For this, we need to construct a group $N$ of type {\bf Am} and a free group $F_g$ fitting into the diagram.
As the cover of $\mathbb{P}^1$ corresponding to $H$ is virtually of Mumford type, the left column is exact with a free group
$F_b$ on $b$ generators and with middle group $N_0=\langle N_{0\bullet}\rangle_{N_0}$. By taking a $*$-tree for
for $N_0$, and making its finite part into a graph of genus $g$ by inserting $g$ edges with trivial stabilisers,
one obtains a group $N$ of type {\bf Am} fitting into the center of (\ref{nineexact}).
Thus, obviously, the middle row is split. 
Further, the vertical map from $F_g$ is the quotient map from $G\to F$,
which implies that the vertical map from $N$ is surjective, as the map $N_0\to H$ is also surjective.
It remains to show that $F_h$ is free, finitely generated.
As $F_c$ is free in $c$ generators, mapping these generators to $F_h$ in such a way that the top right square
commutes, is possible with a map $\alpha\colon F_c\to F_h$. Since the two low rows are split, it follows that $\alpha$ is injective. Thus, $F_h$ is generated by two free groups which intersect trivially, hence  free. 
\end{proof}



If, however, the $G$-cover $\phi\colon X\to Y$ is not of type {\bf Am}, then let $N'$ be an amalgamification of $N$.
The inclusion $N'\subseteq N$ extends to a commutative diagram with exact rows
\begin{align} \label{amalgamifysequence}
\xymatrix{
1 \ar[r]& F_{g'}\ar[d]\ar[r]&N'\ar[d]\ar[r]&G'\ar[d]\ar[r]&1\\
1\ar[r]&F_{g}\ar[r]&N\ar[r]&G\ar[r]&1
}
\end{align}
where $F_{g'}=N'\cap F_g$ is free. Thus, the Galois group $G'$ of the cover of type {\bf Am} can be viewed as a subgroup of $G$. By the construction of amalgamification, $G'$ contains the decomposition groups of $\phi$
and therefore the corresponding cover $\phi'$ contains the ramification datum of $\phi$.
Since the corresponding graph of groups $\mathfrak{G}$ and $\mathfrak{G}'$ are obtained from one another by a finite sequence of 
steps (Am 1) and (Am 2) as in the proof of Lemma \ref{amalgamify}, we can describe, how to obtain $G$ from $G'$. Namely, if $\mathfrak{G}'$ has $s$ loops affected by (Am 1) and $t$ pairs of non-cyclic isomorphic vertex groups $N'_{i_1}\cong N'_{i_2}$ of two different vertices $v_{i_1}$, $v_{i_2}$ obtained by (Am 2), then 
\begin{align}	\label{amalgamgenerate}
G=\langle G', \alpha_1,\dots, \alpha_s, \gamma_{e_1},\dots,\gamma_{e_t}\rangle
\end{align}
subject to relations
\begin{align}	\label{amalgamrelations}
 \gamma_{e_i}N'_{i_1}\gamma_{e_i}^{-1}=N'_{i_2},\quad i=1,\dots,t
\end{align}
(and possibly other relations),
where the $\gamma_{e_i}$ are images of  edges $e_i$ affected by the amalgamification. 

\begin{Prop}
If the Mumford-Hurwitz space $\mathfrak{H}_g(G,{\bf C})$, ${\bf C}=(C_1,\dots,C_r)$, is non-empty, 
then 
$G$ contains a  subgroup
$G'$ containing an element $\gamma_i$ of each $C_i$ such that $\mathfrak{H}_{g'}(G',{\bf C}')$ contains a cover of type {\bf Am} with $*$-graph $\mathfrak{G}'$ (where ${\bf C}'=(C_{G'}(\beta_i))_{i=1,\dots,r+2s+3t}$ are the conjugacy classes of $\gamma_i$ and of conjugates of $\alpha_j$ in $G'$),
and such that $G$ is the quotient of an HNN-extension $N$ of $N'=\pi_1^\BS(\mathfrak{G}^\prime)$
where 
(\ref{amalgamgenerate}) and (\ref{amalgamrelations}) hold with $s+t\le g$. 
\end{Prop}

\begin{proof}
It is clear by the discussion above, that amalgamification gives a cover of type {\bf Am} with the asserted properties. 
\end{proof}

Let now a $G$-cover $\phi$ be given which  allows for a subgroup $G'$ a cover $\phi'$ of type {\bf Am} such that 
(\ref{amalgamgenerate}), (\ref{amalgamrelations}) and $s+t\le g$ hold true.  It is easy to see
$N$ making the diagram (\ref{amalgamifysequence}) commute with exact rows and injective vertical arrows.
By assertion, 
$N$ is a HNN-extension of $N'$ having $N'$ as an amalgamification. But often $N$ is not discretely embeddable into $\PGL_2(K)$.

\section{Realisation of groups}

\subsection{Mumford curves of genus $\ge 2$}
Harbater's result follows easily from Proposition \ref{HM} and can be
improved slightly:

\begin{Thm}     \label{improved}
Every finite group $G$ is the Galois group of a cover
$X\stackrel{G}{\longrightarrow}\mathbb{P}^1$ with a Mumford curve $X$
of genus $g\ge 2$.
\end{Thm}

\begin{proof}
The group $G$, generated by cyclic subgroups $G_1, \dots, G_r$,
is a quotient of the free product $G_1\ast\dots\ast G_r$. This solves, by Proposition
\ref{HM}, the original
inverse Galois problem.
Since, as seen in the proof of Proposition \ref{HM},
 we may realize every cyclic group $C_n$ as a covering group of
a tree with four cusps whose stabilisers are all $C_n$,
the Riemann-Hurwitz-formula tells us
that the genus of the corresponding $C_n$-cover is at least $2$,
if $n>2$. $C_2$-covers are realized as hyperelliptic covers in
a similar way:
$6$ cusps are sufficient for genus $\ge 2$.
\end{proof}

The actual cover can be constructed as in \cite{Harbater1987}.
Translated into the language of $*$-trees, the proof reads:

\smallskip
Let $G$ be generated by two already realised subgroups $H_1$ and
$H_2$ (with Mumford curves of genus $\ge 2$. This means that there
are rational Mumford orbifolds
$\mathcal{S}_1=(\mathbb{P}^1,(\zeta_{11},e_{11}),\dots,(\zeta_{1q},e_{1q}))$ and
$\mathcal{S}_2=(\mathbb{P}^1,(\zeta_{2,q+1},e_{2,q+1}),\dots,(\zeta_{2r},e_{2r}))$,
global charts $\phi_1\colon X_1\stackrel{H_1}{\Pfeil}\mathcal{S}_1$,
$\phi_2\colon X_2\stackrel{H_2}{\Pfeil}\mathcal{S}_2$ with Mumford curves
$X_1$, $X_2$. Let $T_i$ be the $*$-tree for $\phi_i$. The cover
$\phi_i$ induces a morphism of graphs of groups such that the
following diagram is Cartesian:
$$
\xymatrix{
X_i \ar[d]\ar[r]^{H_i} &\mathcal{S}_i\ar[d]\\
\Gamma_i \ar[r]_{H_i} & T_i
}
$$
where the vertical arrows are induced by Tate reduction maps.

\smallskip
Let $T_i^+=T_i\sharf v$ be obtained from $T_i$ by adjoining to any
vertex of $T_i$ an edge $e$ (with trivial group) terminating in $v$
(not belonging to $T_i$). $\Gamma_i^+$ is then defined by the
Cartesian square
$$
\xymatrix{
\Gamma_i^+ \ar[r]\ar[d]_{H_i} &\Gamma_i \ar[d]^{H_i}\\
T_i^+ \ar[r]_{\text{contr.}} &T_i
}
$$
Finally, let
$$
\Delta_i^+:=G/H_i\times\Gamma_i^+,\quad \Delta_i:=G/H_i\times\Gamma_i.
$$
Now, $\Delta_i^+\setminus\Delta_i=G\times e$ for $i=1,2$,
and we can paste $\Delta_i^+\to T_i^+$ along the morphisms
$(\Delta_i^+\setminus\Delta_i\to T_i^+\setminus T_i)
\cong (G\times e\stackrel{\triv}{\Pfeil} e)$, obtaining
$$
\tilde{\phi}\colon \Delta\stackrel{/G}{\Pfeil} T=T_1^+\sharf_eT_2^+
$$
where $T_1^+\sharf_e T_2^+$ means the tree obtained by pasting the trees $T_i^+$ along the edge $e$.
We see that $\Delta$ is connected and $T$ a $*$-tree for an $r$-punctured
orbifold $$\mathcal{S}=\left(\mathbb{P}^1,(\zeta_{ij},e_{ij}),\; i=1,2,\; j=1,\dots,r\right)$$
obtained by
lifting $\tilde{\phi}$ to a chart fitting into
$$
\xymatrix{
X \ar[r]^{G}\ar[d]&\mathcal{S} \ar[d]\\
\Delta \ar[r]_{G} &T
}
$$
with a Mumford curve $X$ and deck group $G$.

\bigskip
The results of
Section \ref{trees-to-graphs} lead to the natural generalisation
of Theorem \ref{improved}
stated in \cite[Corollary 1.2]{vanderPutVoskuil2001} without explicit proof\footnote{which has been added during the time of refereeing the original version of this article \cite[Theorem 1.2]{vanderPutVoskuil2003}}:

\begin{Thm} \label{all-genera}
For every Mumford curve $S$ of genus $g$,
every finite group $G$ is the Galois group of a cover
$X\stackrel{G}{\Pfeil} S$ with a Mumford curve $X$.
\end{Thm}

\begin{proof}
Realise $G$ as a cover of the projective line by a Mumford curve by giving an explicit $*$-tree $T^*$
with fundamental group $N$.
Let $T$ be the finite part of $T^*$.
Inserting $g$ edges between any vertices of the finite part gives a discontinouosly embeddable graph of groups 
with fundamental group\footnote{of type {\bf Am}} $N_g$, thus realising
$G$ over a Mumford curve $S'$ by exhibiting a point of $\mathfrak{M}(N_g)$. 
By forgetting the cover and the $n$ branch points, one sees that there is a 
map 
$$
\mathfrak{M}(N_g)\to \mathfrak{M}_g
$$    
between the ($3g-3+n$)-dimensional moduli space $\mathfrak{M}(N_g)$ and the ($3g-3$)-dimensional space $\mathfrak{M}_g$ of Mumford curves of genus $g$.
This map is easily seen to be surjective.
\end{proof}

\subsection{The full automorphism group}    \label{full}

The extra care taken for realising groups as deck groups in higher genus
allows us to prove

\begin{Thm} \label{fullAut}
Every finite group is isomorphic to the full automorphism group of a Mumford curve.
\end{Thm}

\begin{proof}
Let $G$ be the finite group.
Consider the Mumford-Hurwitz space $\mathfrak{H}_g^h(G,\underline{e})$ of covers $X\stackrel{G}{\Pfeil} Y$
with signature ${\bf e}=(e_1,\dots,e_n)$ with Mumford curves $X$ and $Y$ of
genera $h$ and $g$,
respectively,
where $h$ and $g$ are related by the Riemann-Hurwitz formula
 (the notational redundancy is used here for convenience).
We have a natural map
$$
\Phi\colon\mathcal{H}_g^h(G,{\bf e})\to\mathcal{M}_h,
$$
into the moduli space of all smooth curves of genus $h$,
whose image   $\mathcal{M}_h(G,{\bf e})$ is  of pure dimension $3g-3+n$.
We take a component $\mathcal{M}$ of $\mathcal{M}_h(G,\underline{e})$.
According to \cite[Theorem 5.1]{MSSV2002}, if
$$
(*)\qquad (g,n)\notin\{(2,0), (1,2), (1,1), (0,4), (0,3)\},
$$
then $\mathcal{M}^\circ:=\{X\in \mathcal{M} \mid \Aut X=G\}$ is open and dense in
$\mathcal{M}$. It follows that $\mathcal{M}^\circ$ intersects
$\mathfrak{M}_h(G,{\bf e}):=\Phi(\mathfrak{H}_g^h(G,{\bf e}))$, if and only if $\mathcal{M}$ does,
since $\mathfrak{M}_h(G,{\bf e})$ is of the same pure dimension $3g-3+n$,
provided $\mathfrak{H}_g^h(G,{\bf e})\neq\emptyset$ (Lemma \ref{MumfordHurwitzdimension}). If the
Mumford-Hurwitz space is non-empty, there exist components $\mathcal{M}$
which intersect the Mumford locus $\mathfrak{M}_h(G,{\bf e})$, however. So, we
need only to find for each finite group $G$ a non-empty
Mumford-Hurwitz space $\mathfrak{H}_g^h(G,{\bf e})$ fulfilling $(*)$. But this
is possible by Theorem \ref{all-genera}: just take $g=3$.
\end{proof}

Let us give an
alternative proof using only the realisation over $\mathbb{P}^1$:
Theorem \ref{improved} allows us to realise $G$ as the deck group of a
cover $\phi\colon X\to \mathbb{P}^1$ with a Mumford curve $X$ of genus $g\ge 2$. If $G$ is of order
greater than $5$, Harbater's solution to the inverse Galois problem gives
more than $4$ branch points, and we are done. Otherwise, take the $*$-tree
of $\phi$ and add two edges in such a way that the resulting graph is of genus two.
As it is also a $*$-graph, there is a corresponding Galois cover
$X'\to S'$ fulfilling ($*$),
and we have realised $G$ as $\Aut X'$ with a Mumford curve $X'$ (of genus $\ge 2$).


\subsection*{Acknowledgements}

The author is grateful to Yves Andr\'e for directing the author's
attention toward the theme of this article. He also expresses
thanks to his thesis advisor Frank Herrlich for answering numerous
questions. Thanks also to Harm Voskuil and Fumiharu Kato for fruitful discussions  as well as the unknown referee
for helpful remarks, all of which led to a thoroughly improved rewriting of the original article.
This  would not have been possible
without
Prof.\ Niklaus Kohler  generously letting the author take some time off
for realising this.

\bibliographystyle{gerplain}

\bibliography{retmumf}

\begin{thebibliography}{10}

\bibitem{Andreperiod2003}
{\sc Andre, Yves}: {\em Period mappings and differential equations. From
  $\mathbb{C}$ to $\mathbb{C}_p$ (with appendices by F. Kato and N. Tsuzuki)}.
\newblock Math. Soc. Japan Memoirs, 2003.

\bibitem{Berkovich1990}
{\sc Berkovich, Vladimir~G.}: {\em Spectral Theory and Analytic Geometry over
  Non-Archimedean Fields}.
\newblock Mathematical Surveys and Monographs 33. AMS, 1990.

\bibitem{Berkovich1993}
{\sc Berkovich, Vladimir~G.}: {\em \'Etale cohomology for non-Archimedean
  analytic spaces}.
\newblock Publ. Math. I.H.E.S., 78:5--161, 1993.

\bibitem{Berkovich1999}
{\sc Berkovich, Vladimir~G.}: {\em Smooth $p$-adic analytic spaces are locally
  contractible}.
\newblock Inventiones Mathematicae, 137(1):1--84, 1999.

\bibitem{Bertin1996}
{\sc Bertin, Jos\'e}: {\em Compactifications des sch\'emas de Hurwitz}.
\newblock C. R. Acad. Sci. Paris, L 322(S\'erie 1):1063--1066, 1996.

\bibitem{BradDiss2002}
{\sc Bradley, Patrick~Erik}: {\em $p$-adische Hurwitzr\"aume}.
\newblock Dissertation, Universit\"at Karlsruhe, 2002.

\bibitem{BKVdiscretePGL2}
{\sc {Bradley, P.E.; Kato, F.; Voskuil, H.}}: {\em The Classification of
  $p$-adic Discrete Subgroups of $\rm PGL_2$}.
\newblock In preparation.

\bibitem{deJong1995b}
{\sc {de Jong}, A.~Johan}: {\em \'Etale fundamental groups of non-archimedean
  analytic spaces}.
\newblock Compositio Mathematicae, 97:89--118, 1995.

\bibitem{Fried1993}
{\sc Fried, Michael}: {\em Introduction to MODULAR TOWERS, {\em in} Recent
  Developments in the Inverse Galois Problem}.
\newblock Contemp. Math. 186, 91-171, Seattle, 1993.

\bibitem{Harbater1987}
{\sc Harbater, David}: {\em Galois coverings of the arithmetic line {\em in}
  Number Theory: New York 1984-85}.
\newblock LNM 1240, 165-195, Springer, 1987.

\bibitem{Herrlich1980}
{\sc Herrlich, Frank}: {\em Endlich erzeugbare $p$-adische diskontinuierliche
  Gruppen}.
\newblock Archiv der Mathematik, 35:505--515, 1980.

\bibitem{Herrlich1982}
{\sc Herrlich, Frank}: {\em $p$-adisch discontinuierlich einbettbare Graphen
  von Gruppen}.
\newblock Archiv der Mathematik, 39:204--216, 1982.

\bibitem{Herrlich1984}
{\sc Herrlich~frank}: {\em On the stratification of the moduli space of Mumford
  curves. Groupe d'etude d'Analyse ultrametrique 11e ann\'ee 1983/84, No. 18,
  1-10}.
\newblock 1984.

\bibitem{HerrlichHabil}
{\sc Herrlich, Frank}: {\em Nichtarchimedische Teichm\"ullerr\"aume}.
\newblock Habilitation thesis, Bochum, 1985.

\bibitem{HerrlichHabilArtikel}
{\sc Herrlich, Frank}: {\em Nichtarchimedische Teichm\"ullerr\"aume}.
\newblock Indagationes Mathematicae, 49:149--169, 1987.

\bibitem{Kato2004}
{\sc Kato, Fumiharu}: {\em Non-Archimedean orbifolds covered by Mumford
  curves}.
\newblock J. Algebraic Geometry, 14 (2005):1--34.

\bibitem{Kato99}
{\sc Kato, Fumiharu}: {\em $p$-adic Schwarzian triangle groups of Mumford
  type}.
\newblock Preprint, math.AG/9908174, 1999.

\bibitem{Khramtsov1991}
{\sc Khramtsov, D.~G.}: {\em Finite graphs of groups with isomorphic
  fundamental groups}.
\newblock Algebra Logic, 30(5):389--409, 1991. Translation from Algebra Logika
  30, No. 5, 595-623 (1991).

\bibitem{MSSV2002}
{\sc {K. Magaard, T. Shaska, S. Shpectorov, H. V\"olklein}}: {\em The locus of
  curves with prescribed automorphism group}.
\newblock RIMS, Kyoto Technical Report series, Communications on Arithmetic
  Fundamental Groups and Galois Theory, ed.: H. Nakamura, 2002.

\bibitem{vanderPutVoskuil2001}
{\sc Marius~{van der Put}, Harm~Voskuil}: {\em Mumford coverings of the
  projective line}.
\newblock Preprint, 2001.

\bibitem{vanderPutVoskuil2003}
{\sc Marius~{van der Put}, Harm~Voskuil}: {\em Mumford coverings of the
  projective line}.
\newblock Archiv der Mathematik, 80:98--105, 2003.

\bibitem{WewerDiss}
{\sc Wewers, Stefan}: {\em Construction of Hurwitz spaces}.
\newblock Dissertation, Universit\"at-Gesamthochschule Essen, 1998.

\end{thebibliography}

{\sc Universit\"at Karlsruhe, Institut f\"ur 
Industrielle Bauproduktion,
Eng\-lerstr.~7, D-76128
Karlsruhe, Germany}

e-mail: {\tt bradley@math.uni-karlsruhe.de}


\end{document}